\documentclass [a4 paper,12pt]{article}
\usepackage{mathrsfs}

\usepackage[usenames]{color}
\usepackage{bbm,amsmath,amssymb,amsfonts,amsthm,graphics,graphicx}

\newtheorem{lem}{Lemma}
\newtheorem{thm}[lem]{Theorem}
\newtheorem{cor}[lem]{Corollary}

\newtheorem{pro}{Proposition}

\title{The asymptotic values of\\ the general Zagreb and Randi\'c
indices of\\ trees with bounded maximum degree\footnote{Supported
by NSFC No.10831001, PCSIRT and the ``973" program.}}
\author{\small Xueliang Li, Yiyang Li\\
\small Center for Combinatorics and  LPMC-TJKLC\\
\small Nankai University, Tianjin 300071, China}
\date{ }

\begin{document}

\maketitle

\begin{abstract}
Let $\mathcal {T}^{\Delta}_n$ denote the set of trees of order
$n$, in which the degree of each vertex is bounded by some integer
$\Delta$. Suppose that every tree in $\mathcal {T}^{\Delta}_n$ is
equally likely. We show that the number of vertices of degree $j$
in $\mathcal {T}^{\Delta}_n$ is asymptotically normal with mean
$(\mu_j+o(1))n$ and variance $(\sigma_j+o(1))n$, where $\mu_j$,
$\sigma_j$ are some constants. As a consequence, we give estimate
to the value of the general Zagreb index for almost all trees in
$\mathcal {T}^{\Delta}_n$. Moreover, we obtain that the number of
edges of type $(i,j)$ in $\mathcal {T}^{\Delta}_n$ also has mean
$(\mu_{ij}+o(1))n$ and variance $(\sigma_{ij}+o(1))n$, where an
edge of type $(i,j)$ means that the edge has one end of degree $i$
and the other of degree $j$, and $\mu_{ij}$, $\sigma_{ij}$ are
some constants. Then, we give estimate to the value of the general
Randi\'{c} index for almost all trees in
$\mathcal {T}^{\Delta}_n$.\\[3mm]
\\[3mm]
\end{abstract}

\section{Introduction}

In this paper, we mainly consider trees, in which the degree of
each vertex is bounded by some integer $\Delta$. If $\Delta=1,2$,
the cases are trivial. Thus, we suppose $\Delta\geq 3$ throughout
this paper. Let $\mathcal{T}^{\Delta}_n$ denote the set of trees
with $n$ vertices. We suppose that every tree in
$\mathcal{T}^{\Delta}_n$ is equally likely and $X_n$ is a random
variable, such as the number of vertices of degree $j$, or the
number of edges of type $(i,j)$, each having one end of degree $i$
and the other of degree $j$. It is easy to see that $X_n$ can take
at most $|\mathcal{T}^{\Delta}_n|$ distinct values. We first
introduce two generating functions. Setting
$t_n=|\mathcal{T}^{\Delta}_n|$, we have
$$t(x)=\sum_{n\geq1}t_nx^n,$$
$$t(x,u)=\sum_{n\geq1,k\geq 0}t_{n,k}x^nu^k,$$
where $t_{n,k}$ denotes the number of trees in
$\mathcal{T}^{\Delta}_n$ such that $X_n=k$. Therefore, the
probability of $X_n$ can be defined as
$$\mbox{Pr}[X_n=k]=\frac{t_{n,k}}{t_n}.$$ Note that $t(x,1)=t(x)$.
In \cite{ot}, it is showed that $t_n$ is asymptotically equal to
$\tau\cdot \frac{x_0^{-n}}{n^{5/2}}$, where $\tau$ and $x_0$ are
constants with $x_0\leq1/2$.

In conjunction with the generating functions and asymptotic
analysis, in \cite{dg} and \cite{rs} the authors investigated the
limiting distribution of the number of vertices of given degree
$j$ for trees without degree restriction. By the same method, many
results have been established for other variables, such as the
number of a given {\it path} or {\it  pattern} (see \cite{gk}) for
rooted trees, planar trees, labeled trees {\it et al.} However,
all the statements showed that the limiting distributions are
normal. We refer the readers to \cite{cdkk} and \cite{gk} for
further details.

In this sequel, we follow the method used in \cite{cdkk} and
\cite{dg} to obtain that the distribution of the number of
vertices of degree $j$ for trees in $\mathcal {T}_n^{\Delta}$ is
also asymptotically normal with mean $(\mu_j+o(1))n$ and variance
$(\sigma_j+o(1))n$. Then, we give estimate to the value of the
general Zagreb index for almost all trees in $\mathcal
{T}_n^{\Delta}$. However, for the number of edges of type $(i,j)$,
we only get a weak statement which can not show that the limiting
distribution is normal. Nevertheless, we still can use it to
obtain the asymptotical value of the general Randi\'{c} index for
almost all trees in $\mathcal {T}_n^{\Delta}$.

The definitions of the general Zagreb index and general Randi\'{c}
index will be given in next sections. Many results have been
obtained for the two parameters. We refer the readers to
\cite{lishi} and \cite{lz} for a detailed survey. In this paper we
will show that for the random space $\mathcal {T}_n^{\Delta}$,
each of the indices has a value of $\Theta(n)$ for almost all
trees.

Section $2$ is devoted to a systematic treatment of the number of
vertices of degree $j$ and the general Zagreb index. In Section $3$,
we investigate the number of edges of type $(i,j)$ and the general
Randi\'{c} index.

\section{The number of vertices of degree $j$}

In this section, we first consider the the limiting distribution
of the number of vertices of degree $j$ in
$\mathcal{T}^{\Delta}_n$. Then, as an immediate consequence, we
get the asymptotic value of the general Zagreb index for almost
all trees in $\mathcal{T}^{\Delta}_n$.

In what follows, we introduce some terminology and notation which
will be used in the sequel. For the others not defined here, we
refer to book \cite{hp}.

Analogous to trees, we introduce generating functions for rooted
trees and planted trees. Let $\mathcal{R}^{\Delta}_n$ denote the
set of rooted trees of order $n$ with degrees bounded by an
integer $\Delta$. Setting $r_n=|\mathcal {R}^{\Delta}_n|$, we have
$$r(x)=\sum_{n\geq1}r_n x^n$$ and $$r(x,u)=\sum_{n\geq1, k\geq0}r_{n,k}x^n
u^k,$$ where $r_{n,k}$ denotes the number of trees in
$\mathcal{R}^{\Delta}_n$ such that $X_n$ equals $k$. A {\it
planted tree} is formed by adding a vertex to the root of a rooted
tree. The new vertex is called the {\it plant}, and we never count
it in the sequel. Analogously, let $\mathcal{P}^{\Delta}_n$ denote
the set of planted trees with $n$ vertices of bounded maximum
degree $\Delta$. Setting $p_n=|\mathcal{P}^{\Delta}_n|$, we have
$$p(x)=\sum_{n\geq1}p_n x^n$$ and $$p(x,u)=\sum_{n\geq 1,
k\geq0}p_{n,k}x^n u^k,$$ where $p_{n,k}$ denotes the number of
trees in $\mathcal{P}^{\Delta}_n$ such that $X_n$ equals $k$. By
the definition of planted trees, one can readily see that
$p(x,1)=p(x)=r(x,1)=r(x)$.

Furthermore, we introduce another generating function
$p^{(\Delta-1)}(x)$. Denote $p_n^{(\Delta-1)}$ as the number of
planted trees such that the degree of the root is bounded by
$\Delta-1$, while the degrees of other vertices are bounded by
$\Delta$. Then, we define
$$p^{(\Delta-1)}(x)=\sum_{n\geq1}p_n^{(\Delta-1)}x^n.$$

In \cite{ot}, Otter showed that there exists a number $x_0$ such
that \begin{equation}\label{expand}p(x)=b_1+b_2\sqrt{x_0-x}+b_3
(x_0-x)+\cdots,\end{equation} where $b_1,b_2,b_3$ are some
constants not equal to zero. Evidently, $p(x_0)=b_1$ and for any
$|x|\leq x_0$, $p(x)$ is convergent. For any $\Delta$,
$x_0\leq1/2$; particularly, if $\Delta=4$, $x_0\approx 0.3551817$
and $p(x_0)\approx1.117421$. Moreover,
\begin{equation}\label{De-1}
p^{(\Delta-1)}(x_0)=1.
\end{equation}
We refer the readers to \cite{ot} for more details.

The proofs of our main results in this paper ultimately rely on
the following lemma, due to Chyzak {\it et al.} \cite{cdkk} and
Drmota \cite{d}. We first introduce some notation.

Let $\textbf{y}(x,u)=(y_1(x,u),\ldots,y_N(x,u))^T$ be a column
vector. We suppose that $G(x,\textbf{y},u)$ is an analytic function
with non-negative integer Taylor coefficients. $G(x,\textbf{y},u)$
can be expanded as
$$ G(x,\textbf{y},u)=\sum_{n\geq1,k\geq0}g_{n,k}x^nu^k.$$
Let $X_n$ denote a random variable with probability
\begin{equation}\label{rv}
\mbox{Pr}[X_n=k]=\frac{g_{n,k}}{g_n},
\end{equation}
where $g_n=\sum_kg_{n,k}$.
\begin{lem}\label{mainlem}
Let
$\textbf{F}(x,\textbf{y},u)=(F_1(x,\textbf{y},u),\ldots,F_N(x,\textbf{y},u))^T$
be functions analytic around $x=0$,
$\textbf{y}=(y_1,\ldots,y_N)^T=\textbf{0}$, $u=0$, with Taylor
coefficients all are non-negative integers. Suppose
$\textbf{F}(0,\textbf{y},u)=\textbf{0}$,
$\textbf{F}(x,\textbf{0},u)\neq\textbf{0}$,
$\textbf{F}_x(x,\textbf{y},u)\neq \textbf{0}$, and for some $j$,
$\textbf{F}_{y_jy_j}(x,\textbf{y},u)\neq\textbf{0}$. Furthermore,
assume that $x=x_0$ together with $\textbf{y}=\textbf{y}_0$ is a
non-negative solution of the system of equations
\begin{align}
&\textbf{y}=\textbf{F}(x,\textbf{y},1)\\
&0=\mbox{det}(\textbf{I}-\textbf{F}_\textbf{y}(x,\textbf{y},1))\label{det}
\end{align}
inside the region of convergence of $\textbf{F}$, $\textbf{I}$ is
the unit matrix. Let $\textbf{y}=(y_1(x,u),\ldots,y_N(x,u))^T$
denote the analytic solution of the system
\begin{equation}\label{fs}
\textbf{y}=\textbf{F}(x,\textbf{y},u)
\end{equation}
with $\textbf{y}(0,u)=\textbf{0}$.

If the dependency graph $G_{\textbf{F}}$ of the function system
Equ.(\ref{fs}) is strongly connected, then there exist functions
$f(u)$ and $g_i(x,u)$, $h_i(x,u)$ ($1\leq i\leq N$) which are
analytic around $x=x_0$, $u=1$, such that
\begin{equation}\label{sqrt}
y_i(x,u)=g_i(x,u)-h_i(x,u)\sqrt{1-\frac{x}{f(u)}}
\end{equation}
is analytically continued around $u=1$, $x=f(u)$ with
$\mbox{arg}(x-f(u))\neq 0$, where $x=f(u)$ together with
$y=y(f(u),u)$ is the solution of the extended system
\begin{align}
&\textbf{y}=\textbf{F}(x,\textbf{y},u)\\
&0=\mbox{det}(\textbf{I}-\textbf{F}_\textbf{y}(x,\textbf{y},u))\label{det2}.
\end{align}

Moreover, let $G(x,\textbf{y},u)$ be an analytic function with
non-negative Taylor coefficients such that the point
$(x_0,\textbf{y}(x_0,1),1)$ is contained in the region of
convergence. Finally, let $X_n$ be the random variable defined in
Equ.(\ref{rv}). Then the random variable $X_n$ is asymptotically
normal with mean
$$E(X_n)=\mu n+O(1)\mbox{  } (n\rightarrow \infty),$$ and variance
$$Var(X_n)=\sigma n+O(1)\mbox{ } (n\rightarrow \infty)$$
with $\mu=\frac{-f'(1)}{f(1)}$.
\end{lem}

\noindent {\bf Remark 1:} We say that the {\it dependency graph}
$G_{\textbf{F}}$ of $\textbf{y}=\textbf{F}(x,\textbf{y},u)$ is
strongly connected if there is no subsystem of equations that can
be solved independently from others. If $G_\textbf{F}$ is strongly
connected, then
$\textbf{I}-\textbf{F}_\textbf{y}(x_0,\textbf{y}_0,1)$ has rank
$N-1$. Suppose that $\mathbf{v}^T$ is a vector with
$\mathbf{v}^T(\mathbf{I}-\mathbf{F}_\mathbf{y}(x_0,\mathbf{y_0},1))=0$.
Then, $\mu=\frac{\mathbf{v}^T(\mathbf{F}_u(x_0,\mathbf{y_0},1))}
{x_0\mathbf{v}^T(\mathbf{F}_x(x_0,\mathbf{y_0},1))}$. We refer the
readers to \cite{cdkk, d} for more details.\\

In what follows, we shall use the above lemma to investigate the
number of vertices of degree $j$ in $T^{\Delta}_n$, where $j$ is a
given integer.

Firstly, we focus on the planted trees. There appears an
expression of the form $Z(S_n, f(x,u))$ (or $f(x)$), which is the
substitution of the counting series $f(x,u)$ (or $f(x)$) into the
cycle index $Z(S_n)$ of the symmetric group $S_n$. This involves
replacing each variable $s_i$ in $Z(S_n)$ by $f(x^i,u^i)$ (or
$f(x^i)$). For instance, if $n=3$, then
$Z(S_3)=(1/3!)(s_1^3+3s_1s_2+2s_3)$, and
$Z(S_3,f(x,u))=(1/3!)(f(x,u)^3+3f(x,u)f(x^2,u^2)+2f(x^3,u^3))$. We
refer the readers to \cite{hp} for details.

Note that a planted tree with a root of degree $k$ can be viewed
as a root vertex attached by $k-1$ planted trees. Employing the
classic P\'{o}lya enumeration theorem, we have $Z(S_{k-1};p(x))$
as the counting series of the planted trees whose roots have
degree $k$, and the coefficient of $x^p$ in $x\cdot
Z(S_{k-1};p(x))$ is the number of planted trees with $p$ vertices
(see \cite{hp} p.51--54). Therefore,
$$p(x)=x\cdot\sum_{k=0}^{\Delta-1}Z(S_k;p(x)),$$ and
$$p^{(\Delta-1)}(x)=x\cdot\sum_{k=0}^{\Delta-2}Z(S_k;p(x)).$$

By the same method, we can obtain that
\begin{equation}\label{pj}
p(x,u)=x\cdot\sum_{l=1}^{\Delta}Z(S_{l-1};
p(x,u))+x(u-1)Z(S_{j-1};p(x,u)),
\end{equation}
where the last term $(xu-x)Z(S_{j-1};p(x,u))$ serves to count the
vertices of degree $j$ when the root of a planted tree is of
degree $j$. Then, we show that this equation satisfies the
conditions of Lemma \ref{mainlem}. Suppose $p(x,u)=F(x,p(x,u),u)$.
It is well-known that the partial derivative of $Z(S_n;\cdot)$
enjoys (see \cite{dg})
\begin{equation}\label{pdz}
\frac{\partial}{\partial s_1}Z(S_n;
s_1,\ldots,s_n)=Z(S_{n-1};s_1,\ldots,s_{n-1}).
\end{equation}
One can readily see that
$$F_p(x_0,p(x_0,1),1)=x_0\sum_{k=0}^{\Delta-2}Z(S_k; p(x_0,1))=p^{(\Delta-1)}(x_0)=1.$$
The other conditions are easy to be illustrated. Thus we have that
$p(x,u)$ is in the form of
\begin{equation}\label{eqp}
p(x,u)=g_1(x,u)-h_1(x,u)\sqrt{1-\frac{x}{f(u)}},
\end{equation}
where $g_1(x,u)$, $h_1(x,u)$ and $f(u)$ are analytic around $x=x_0$
and $u=1$, and $p(x,u)$ is analytically continued around $u=1$,
$x=f(u)$ with $\mbox{arg}(x-f(u))\neq 0$. From Equ.(\ref{expand}),
we can see $f(1)=x_0$.

Analogous to Equ.(\ref{pj}), for rooted trees, it follows that
\begin{equation}\label{rp}
r(x,u)=x\cdot\sum_{l\geq0}^{\Delta}Z(S_{l};p(x,u))+(xu-x)Z(S_{j};p(x,u)).
\end{equation}

For trees, however, in \cite{ot} the author obtained that
\begin{eqnarray}\label{trp}
t(x)=r(x)-\frac{1}{2}p(x)^2+\frac{1}{2}p(x^2),
\end{eqnarray}
and $t(x)$ can be expanded as
\begin{equation}\label{expandt}
t(x)=c_1+c_2(1-\frac{x}{x_0})+c_3{(1-\frac{x}{x_0})}^{3}+\cdots,
\end{equation}
where $c_1, c_2, c_3$ are some constants not equal to zero.
Furthermore, we can also obtain a similar equation for $t(x,u)$.
We introduce a useful lemma due to Otter \cite{ot}.

Two edges in a tree are {\it similar}, if they are the same under
some automorphism of the tree. To {\it join} two planted trees is
to connect the two roots with a new edge and get rid of the two
plants. If the two panted trees are the same, we say that the new
edge is {\it symmetric}.
\begin{lem}\label{otter}
For any tree, the number of rooted trees corresponding to this
tree minus the number of nonsimilar edges (except for the
symmetric edge) is the number $1$.
\end{lem}

Note that, if we delete any one edge from a similar set in a tree,
the yielded trees are the same two trees. Hence, different pairs
of planted trees correspond to nonsimilar edges. Now, we have
\begin{equation}\label{t(x,u)}
t(x,u)=r(x,u)-\frac{1}{2}p(x,u)^2+\frac{1}{2}p(x^2,u^2).
\end{equation}
Here, we should notice that $t(x,u)$ is a function in $p(x,u)$, but
its Taylor coefficient of $p(x,u)$ is not sure to be non-negative.
Thus, we could not use Lemma \ref{mainlem} for $t(x,u)$.

However, Equ.(\ref{t(x,u)}) together with Equs.(\ref{eqp}) and
(\ref{rp}) gives
$$t(x,u)=\widetilde{g}(x,u)-\widetilde{h}(x,u)\sqrt{1-\frac{x}{f(u)}},$$
where $\widetilde{g}(x,u)$, $\widetilde{h}(x,u)$ are analytic
around $x=f(1)$ and $u=1$, and $t(x,u)$ is also analytically
continued around $x=f(u)$ and $u=1$ with $\mbox{arg}(x-f(u))\neq
0$.

In what follows, we shall show that $\widetilde{h}(f(u),u)=0$ around
$u=1$. We have
\begin{align*}
t(x,u)&=p(x,u)-x(u-1)Z(S_{j-1};p(x,u))+x\cdot
Z(S_{\Delta};p(x,u))\\
&+x(u-1)Z(S_j;
p(x,u))-\frac{1}{2}p(x,u)^2+\frac{1}{2}p(x^2,u^2)\\
&=g_1-h_1\sqrt{1-\frac{x}{f(u)}}-x(u-1)Z(S_{j-1};g_1-h_1\sqrt{1-\frac{x}{f(u)}})\\
&+x\cdot Z(S_{\Delta};g_1-h_1\sqrt{1-\frac{x}{f(u)}}) +x(u-1)Z(S_j;
g_1-h_1\sqrt{1-\frac{x}{f(u)}})\\
&-\frac{1}{2}(g_1^2+h_1^2(1-\frac{x}{f(u)}))\\
&+g_1h_1\sqrt{1-\frac{x}{f(u)}}+p(x^2,u^2).
\end{align*}
Then, by means of Taylor's theorem we get
$$Z(S_k;g_1-h_1\sqrt{1-\frac{x}{f(u)}})=
\sum_{i=0}^kZ^{(i)}(S_k;g_1)h_1^i(1-x/f(u))^{i/2}\frac{(-1)^i}{i!},
$$ where $Z^{(i)}$ denotes the $i$-th derivative with respect to the
cycle index $s_1$. Thus, we obtain
$$\widetilde{h}=h_1(1-g_1+x\cdot
Z(S_{\Delta-1};g_1)+x(u-1)Z(S_{j-1};g_1)-x(u-1)Z(S_{j-2};g_1)).$$

On the other hand,  from Lemma \ref{mainlem}, note that $x=f(u)$ and
$p(x,u)=g_1(f(u),u)$ are the solutions of
\begin{align*}
p&=x\cdot\sum_{l=1}^{\Delta}Z(S_{l-1}; p)+x(u-1)Z(S_{j-1};p),\\
1&=x\cdot\sum_{l=1}^{\Delta-1}Z(S_{l-1};p)+x(u-1)Z'(S_{j-1};p),
\end{align*}
which yields $$g_1(f(u),u)=1+f(u)\cdot
Z(S_{\Delta-1};g_1)+f(u)(u-1)(Z(S_{j-1};g_1)-Z(S_{j-2};g_1)),$$
that is, $$\widetilde{h}(f(u),u)=0.$$  Then by setting
$h(x,u)=\widetilde{h}(x,u)(1-\frac{x}{f(u)})^{-1}$, for some
$g(x,u)$ it follows that
$$t(x,u)=g(x,u)-h(x,u)(1-\frac{x}{f(u)})^{3/2},$$ around $x=f(u)$
and $u=1$ with $\mbox{arg}(x-f(u))\neq 0$. Moreover, we can see
$h(x,u)$ is analytic around $x=f(1)$ and $u=1$. By
Equ.(\ref{expandt}), we have $h(f(1),1)\neq 0$, and thus
$h(f(u),u)\neq 0$ around $u=1$.

Next, we need the following proposition from \cite{dg}, which can
be proved by a transfer lemma of Flajolet and Odlyzko \cite{fo}
and Cauchy's formula. We refer the readers to \cite{d, dg} and
\cite{fo} for more details.

\begin{pro}
Suppose $y(x,u)=\sum y_{nm}x^nu^m$ is an analytic function with
$y_{nm}\geq 0$. There exist functions $g(x,u)$, $h(x,u)$ and
$f(u)$ which are analytic around $x=x_0=f(1)$ and $u=1$, $x_0$ is
the radius of convergence of $y(x,1)$, $y(x,u)$ is analytically
continued in the region $|x-f(u)|<\eta$, $arg(x-f(u))\neq 0$ and
$|u-1|<\eta$, where $\eta$ is sufficiently small, and
$$y(x,u)=g(x,u)-h(x,u){(1-\frac{x}{f(u)})^{3/2}}.$$
Then $y_n(u)=\sum_my_{nm}u^m$ is asymptotically given by
$$y_n(u)=\frac{3h(f(u),u)}{4\sqrt{\pi}n^{5/2}}f(u)^{-n+1}+O(\frac{f(u)^{-n}}{n^{7/2}})
$$ uniformly for $|u-1|<\eta$.
If $h(f(1),1)\neq 0$ and $X_n$ is defined as Equ.(\ref{rv}) for
$y(x,u)$, then $X_n$ is asymptotically normal with mean
$(\mu+o(1))n$ and variance $(\sigma+o(1))n$ where $\mu$ and
$\sigma$ are some constants.
\end{pro}

We can see that all the conditions hold for $t(x,u)$. Then, for
the number of vertices of degree $j$, the following result is
immediate.
\begin{thm}\label{thm1}
Suppose $j$ is an integer. Let $X_n$ be the number of vertices of
degree $j$ in $\mathcal{T}_n$. Then, $X_n$ is asymptotically
normally distributed with mean value $(\mu_j+o(1)) n$ and variance
$(\sigma_j+o(1)) n$, where $\mu_j$ and $\sigma_j$ are some constants
to every $j$.
\end{thm}

Following book \cite{BB}, we will say that {\it almost every}
(a.e.) graph in a random graph space $\mathcal{G}_n$ has a certain
property $Q$ if the probability $\mbox{Pr}(Q)$ in $\mathcal{G}_n$
converges to 1 as $n$ tends to infinity. Occasionally, we shall
write {\it almost all} instead of almost every.

For the number of vertices of degree $j$, by Chebyshev inequality
one can get that
$$ \mbox{Pr}\big[\big|X_n-E(X_n)\big|>n^{3/4}\big]\leq \frac{Var
X_n}{n^{3/2}}\rightarrow 0 \mbox{ as } n\rightarrow \infty.$$
Therefore, $X_n=(\mu_j+o(1)) n \mbox{  a.e.} $ Then, an immediate
consequence is the following.

\begin{cor}\label{ch}
For almost all trees in $\mathcal{T}^{\Delta}_n$, the number of
vertex of degree $j$ is $(\mu_j+o(1)) n$.
\end{cor}

Now, we discuss the general Zagreb index. Let $G=(\mathcal
{V},\mathcal {E})$ be a graph with vertex set $\mathcal {V}$ and
edge set $\mathcal{E}$. The {\it general Zagreb index} was
introduced by Li {\it et al.} \cite{lz}, where they call it {\it
the zeroth order general Randi\'c index}, and is defined to be the
sum of powers of degree, i.e.,
$$D_{\alpha}=\sum_{u\in \mathcal {V}}{d_u^{\alpha}},$$ where $\alpha$ is some
real number, $d_u$ is the degree of vertex $u$. Many results have
been obtained for this variable. Particularly, if $\alpha=-1$,
$D_{-1}$ is called the {\it inverse degree}, and if $\alpha=2$,
$D_2$ is known as the {\it first Zagreb index} \cite{gt}.

For a tree with $n$ vertices, if $\alpha=0$ then $D_0=n$. So, we
focus on the case $\alpha\neq 0$, and establish estimate to the
value of $D_{\alpha}$ for almost all trees in $\mathcal
{T}^{\Delta}_n$.

Since the degrees of the tree in $\mathcal{T}_n^{\Delta}$ are
bounded by $\Delta$, we can obtain that, for almost all trees,
$$D_{\alpha}=\sum_{j=1}^{\Delta} j^{\alpha}\cdot (\mu_j+o(1))n.$$
For convenience, set $d_{\alpha}=\sum_{j=1}^{\Delta}
j^{\alpha}\cdot \mu_j.$

\begin{cor}
For almost all trees in $\mathcal{T}_n^{\Delta}$, the value of the
general Zagreb index enjoys
$$D_{\alpha}=(d_{\alpha}+o(1))n,$$ where $d_{\alpha}$ is a constant.
\end{cor}

\section{The number of edges of type $(i,j)$}

We start this section by counting the number of edges of type
$(i,j)$. Without loss of generality, suppose $i\leq j$. Since
there is only one tree with an edge of type $(1,1)$, we always
assume $j>1$.

In this section, we still use the same notation as in Section $2$.
We also use $X_n$ to denote the number of edges of type $(i,j)$,
which would not make any ambiguity. Split up $P_n^{\Delta}$ into
$\Delta$ subsets according to the degrees of the roots, and let
$a_k(x,u)$ (or $a_k$) be the generating function corresponding to
each subset. Then, we have $a_1(x,u)=x$ and
$$x+a_2(x,u)+\cdots+a_\Delta(x,u)=p(x,u).$$

Firstly, we establish the functions system and use Lemma
\ref{mainlem} to get equations for $a_k(x,u)$ and $p(x,u)$ in the
form of Equ.(\ref{sqrt}). Analogous to Equ.(\ref{pj}), we can
obtain equations below
\begin{align}\label{fs4}
a_2(x,u)=&x\cdot p(x,u),\nonumber\\
\cdots&\cdots\nonumber\\
a_i(x,u)=&x\cdot \sum_{\ell_1+\ell_2=i-1}Z(S_{\ell_1}; p(x,u)-a_j)\cdot Z(S_{\ell_2};a_j)u^{\ell_2},\nonumber\\
\cdots&\cdots\nonumber\\
a_j(x,u)=&x\cdot \sum_{m_1+m_2=j-1}Z(S_{m_1}; p(x,u)-a_i)\cdot Z(S_{m_2};a_i)u^{m_2},\nonumber\\
\cdots&\cdots\nonumber\\
a_\Delta(x,u)=&x\cdot Z(S_{\Delta-1};p(x,u)).\nonumber
\end{align}
With respect to $a_i(x,u)$, the root of a tree is attached by
$i-1$ planted trees. Then, suppose there are $\ell_2$ planted
trees attached to the root that has degree $j$ while the other
$\ell_1$ planted trees are not, for which, in the above equations
the term $Z(S_{\ell_1}; p(x,u)-a_j)\cdot
Z(S_{\ell_2};a_j)u^{\ell_2}$ is to treat this case. Hence, the
above functions system follows.

Moreover, we shall show that these functions satisfy all the
conditions of Lemma \ref{mainlem}. Since the others are easy to
verify, we only show that Equ.(\ref{det}) holds. Set
$\mathbf{a}(x,u)=(a_2,\ldots,a_{\Delta})$. In fact, by using
Equ.(\ref{pdz}), for $k\geq 2$ it is easy to get
\begin{align*}
\mathbf{F}_{a_k}(x_0,\mathbf{a}(x_0,1),1)=\begin{pmatrix} x_0\\
\cdots\\
\\x_0\cdot Z(S_{j-2};p(x_0,1))\\
\cdots
\\x_0\cdot Z(S_{\Delta-2};p(x_0,1))
\end{pmatrix}.
\end{align*}
Recalling that $p^{(\Delta-1)}(x_0,1)=p^{(\Delta-1)}(x_0)=1$, one
can readily see that
$$\mbox{det}(\mathbf{I}-\mathbf{F}_\mathbf{a}(x_0,\mathbf{a}(x_0,1),1))=0.$$
Moreover, we have $\mathbf{v}^T=(1,\ldots,1)$. Then from Lemma
\ref{mainlem}, we have\begin{equation}\label{ag_k}
a_k(x,u)=g_k(x,u)-h_k(x,u)\sqrt{1-\frac{x}{f(u)}}
\end{equation}
where $2\leq k\leq \Delta$, and thus
\begin{equation}\label{pg_0}
p(x,u)=g_0(x,u)-h_0(x,u)\sqrt{1-\frac{x}{f(u)}},
\end{equation}
where $g_k$, $h_k$, $g_0$, $h_0$ and $f(u)$ are as required in
Lemma \ref{mainlem}. Since $p(x,u)$ is the sum of $a_k$'s, which
is the function in $\mathbf{a}(x,u)$, $x$ and $u$, by Lemma
\ref{mainlem} we have
$\frac{-f'(1)}{f(1)}=\frac{1}{x_0}\frac{\mathbf{v}^T(\mathbf{F}_u(x_0,\mathbf{a}_0,1))}
{\mathbf{v}^T(\mathbf{F}_x(x_0,\mathbf{a}_0,1))}$ and $f(1)=x_0$.
From the functions system of $a_k$'s, we get that
$\mathbf{F}_u(x_0,\mathbf{a}_0,1)$ and
$\mathbf{F}_{x}(x_0,\mathbf{a}_0,1)$ are positive and thus
$\frac{-f'(1)}{f(1)}>0$.

For the rooted trees, we have
\begin{align*}
r(x,u)&=x\cdot
\sum_{k=1}^{\Delta}Z(S_{k};p(x,u))+x\sum_{l_1+l_2=i}Z(S_{l_1};
p(x,u)-a_j)\cdot Z(S_{l_2};a_j)(u^{l_2}-1)\\
&+x\cdot \sum_{m_1+m_2=j}Z(S_{m_1}; p(x,u)-a_i)\cdot
Z(S_{m_2};a_i)(u^{m_2}-1).
\end{align*}

Note that if we join two planted trees with roots of degree $i$
and $j$, respectively, then the number of edges of type $(i,j)$ in
the new tree is counted by $a_i(x,u)a_j(x,u)u$. Therefore, it
follows that if $i<j$, then
\begin{align*}
t(x,u)=&r(x,u)-\frac{1}{2}p(x,u)^2+\frac{1}{2}p(x^2,u^2)+a_ia_j(1-u),\end{align*}
and if $i=j$, then
\begin{align*}
t(x,u)=&r(x,u)-\frac{1}{2}p(x,u)^2+\frac{1}{2}p(x^2,u^2)
+\frac{1}{2}a_i^2(1-u)-\frac{1}{2}a_i(x^2,u^2)(1-u).
\end{align*}
By using Equs.(\ref{pg_0}) and (\ref{ag_k}), we always obtain that
$t(x,u)$ is in the form of
\begin{align*}
t(x,u)&=\overline{g}-\overline{h}\sqrt{1-\frac{x}{f(u)}}.
\end{align*}
Surely, we can proceed to show that $\overline{h}(f(u),u)=0$
around $u=1$ as previous and obtain that the random variable $X_n$
to the edges of type $(i,j)$ is still asymptotically normal. But
it is much involved in this case. Thus, we introduce the following
lemma (see \cite{gk}), which can give us a weak result.
\begin{lem}\label{lem3}
Suppose $y(x,u)$ has the form
$$y(x,u)=g(x,u)-h(x,u)\sqrt{1-\frac{x}{f(u)}}$$
where $g(x,u)$, $h(x,u)$ and $f(u)$ are analytic functions around
$x=f(1)$ and $u=1$ that satisfy $h(f(1),1)=0$, $h_x(f(1),1)\neq
0$, $f(1)>0$ and $f'(1)<0$. Furthermore, $x=f(u)$ is the only
singularity on the cycle $|x|=|f(u)|$ for $u$ is close to $1$.
Suppose $X_n$ is defined as Equ.(\ref{rv}) to $y(x,u)$. Then,
$E(X_n)=(\mu+o(1))n$ and $Var(X_n)=(\sigma+o(1))n$, where
$\mu=-f'(1)/f(1)$ and $\sigma=\mu^2+\mu-f''(1)/f(1)$.
\end{lem}

\noindent {\bf Remark 2:} This result does not tell us that the
limiting distribution is asymptotically normal. If $h(f(1),1)\neq
0$, this lemma is trivial by Lemma \ref{mainlem}, and if
$h(f(u),u)=0$, we can still get that the limiting distribution is normal.\\

For $t(x,u)$, since $t(x,1)=t(x)$, we have that
$\overline{h}(f(1),1)=0$ and $\overline{h}_x(f(1),1)\neq0$.
Moreover, the other conditions in Lemma \ref{lem3} immediately
follow from Equs.(\ref{ag_k}) and (\ref{pg_0}). Then, we can
establish the following theorem.
\begin{thm}
Let $X_n$ be the number of edges of type $(i,j)$ in
$\mathcal{T}_n^{\Delta}$. Then,
$$E(X_n)=(\mu_{ij}+o(1)) n$$ and
$$Var(X_n)=(\sigma_{ij}+o(1))n$$ where $\mu_{ij}$ and
$\sigma_{ij}$ are some constants to every type $(i,j)$.
\end{thm}

Consequently, by Chebyshev inequality we have the following
result.
\begin{cor}\label{corlast}
For almost all trees in $\mathcal{T}_n^{\Delta}$, the number of
edges of type $(i,j)$ equals $(\mu_{ij}+o(1))n$.
\end{cor}

Now, we can give estimate to the value of the general Randi\'{c}
index for trees in $\mathcal {T}^{\Delta}_n$. The {\it general
Randi\'{c} index} is defined as
$$R_{\beta}=\sum_{uv\in \mathcal {E}} (d_ud_v)^{\beta},$$
where $d_u$, $d_v$ are the degree of $u$, $v$, respectively. If
$\beta=-1/2$, the index is called the classic {\it Randi\'{c}
index} \cite{r}. If $\beta=1$, $R_{1}$ is known as the {\it second
Zagreb index} \cite{gt}. We refer the readers to \cite{lishi} for
a detailed survey. Moreover, for a tree with $n$ vertices, if
$\beta=0$ then $R_0=n-1$. Thus, we suppose $\beta\neq 0$. We shall
get the estimate of $R_{\beta}$ for almost all trees.

By Corollary \ref{corlast}, for trees in $\mathcal {T}^{\Delta}_n$,
we can obtain that
$$R_{\beta}=\sum_{i\leq j\leq \Delta}(ij)^{\beta}\cdot(\mu_{ij}+o(1))n \mbox{  a.e.}$$
Denote $\sum_{i\leq j\leq \Delta}(ij)^{\beta}\cdot\mu_{ij}$ by
$r_\beta$. Then the following result is immediate.

\begin{cor}
For almost all trees in $\mathcal{T}^{\Delta}_n$, the general
Randi\'{c} index $R_{\beta}$ equals $(r_\beta+o(1))n$, where
$r_\beta$ is a constant.
\end{cor}

\section {Concluding remark}

Although the general Zagreb index and the general Randi\'c index for
trees and general graphs have been studied extensively, there are
very few results from the asymptotic point of view to study them.
Almost all known results are about extremal values of the indices
and trees or graphs that attain the extremal values. It is known
(see \cite{lishi}) that for trees the maximum of the first Zagreb
index and the Randi\'c index is, respectively, $n(n-1)$ (attained by
the star $S_n$) and $\frac {n-3} {2} +\sqrt {2}$ (attained by the
path $P_n$), while the minimum of them is, respectively, $4n-6$
(attained by the path $P_n$) and $\sqrt {n-1}$ (attained by the star
$S_n$). Our results show that for almost all trees of bounded
maximum degree the values of the indices are linear in the order $n$
of the trees.

\end{document}